# An optimal control-based numerical method for scalar transmission problems with sign-changing coefficients

Patrick Ciarlet JR[1], David Lassounon[2], Mahran Rihani[3]

[1] POEMS, Ensta Paris, Institut Polytechnique de Paris, 91128 Palaiseau , France.

[2] Institut de Recherche Mathématique de Rennes (IRMAR), Université Européenne de Bretagne, 20 avenue des Buttes de Coësmes, CS 70839, 35708 Rennes Cédex 7 , France.

[3] CMAP, Ecole Polytechnique, Institut Polytechnique de Paris, 91128 Palaiseau, France.

E-mails: patrick.ciarlet@ensta.fr, enagnon-david.lassounon@insa-rennes.fr, mahran.rihani@polytechnique.edu.

(May 12, 2022)

**Abstract:** In this work, we present a new numerical method for solving the scalar transmission problem with sign-changing coefficients. In electromagnetism, such a transmission problem can occur if the domain of interest is made of a classical dielectric material and a metal or a metamaterial, with for instance an electric permittivity that is strictly negative in the metal or metamaterial. The method is based on an optimal control reformulation of the problem. Contrary to other existing approaches, the convergence of this method is proved without any restrictive condition. In particular, no condition is imposed on the a priori regularity of the solution to the problem, and no condition is imposed on the meshes, other than that they fit with the interface between the two media. Our results are illustrated by some (2D) numerical experiments.

**Key words:** transmission problem, sign-changing coefficients, fictitious domain methods, optimal control.

## 1 Introduction

In the present paper, we study the numerical approximation of the scalar transmission problem with sign-changing coefficients in $\mathbb{R}^d$, for $d \in \{2,3\}$. To fix ideas, let $\Omega$ be an open, bounded, connected subset of $\mathbb{R}^d$ with a Lipschitz boundary, in other words a *domain* of $\mathbb{R}^d$. Further, consider that $\Omega$ is equal to the union of two disjoint (sub)domains $\Omega_1, \Omega_2$. We denote the interface by $\Sigma = \partial\Omega_1 \cap \partial\Omega_2$ (see Figure 1 for an example), and we assume that $\text{meas}_{\partial\Omega}(\partial\Omega_2 \setminus \Sigma) > 0$.

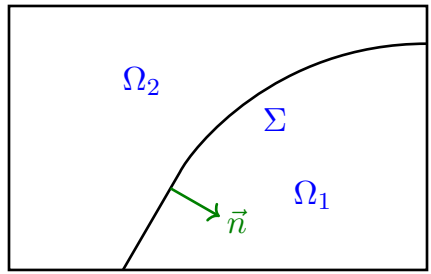

Figure 1: Example of geometry.

We also introduce a coefficient $\varepsilon \in L^\infty(\Omega)$ such that $\varepsilon_1 = \varepsilon_{|\Omega_1} \geq \varepsilon_+ > 0$ a.e. in $\Omega_1$ and $\varepsilon_2 = \varepsilon_{|\Omega_2} \leq \varepsilon_- < 0$ in a.e. in $\Omega_2$. Here $\varepsilon_+$ and $\varepsilon_-$ are two real constants. It will be useful to introduce the contrasts $\kappa_\varepsilon^1 := \varepsilon_1^- / \varepsilon_2^+$ and $\kappa_\varepsilon^2 := \varepsilon_2^- / \varepsilon_1^+$ where $\varepsilon_1^\pm$ and $\varepsilon_2^\pm$ are defined as follows:

$$\varepsilon_1^+ := \sup_{\Omega_1} \varepsilon_1, \quad \varepsilon_1^- := \inf_{\Omega_1} \varepsilon_1, \quad \varepsilon_2^+ := \sup_{\Omega_2} |\varepsilon_2| \text{ and } \varepsilon_2^- := \inf_{\Omega_2} |\varepsilon_2|.$$

Note that in the particular case where $\varepsilon$ is piecewise constant, we have $\kappa_\varepsilon^1 = 1/\kappa_\varepsilon^2$.

**Remark 1.1.** *In principle, $\varepsilon$ could be a symmetric tensor-valued coefficient, ie. $\varepsilon = (\varepsilon_{ij})_{1 \leq i,j \leq d}$ with $\varepsilon_{ij} \in L^\infty(\Omega)$ for all $1 \leq i, j \leq d$, and such that*

$$\exists \varepsilon_+ > 0, \ \forall \boldsymbol{z} \in \mathbb{R}^d, \quad \varepsilon_+ |\boldsymbol{z}|^2 \leq \varepsilon \boldsymbol{z} \cdot \boldsymbol{z} \ a.e. \ in \ \Omega_1 \,;$$
$$\exists \varepsilon_- > 0, \ \forall \boldsymbol{z} \in \mathbb{R}^d, \quad \varepsilon_- |\boldsymbol{z}|^2 \leq -\varepsilon \boldsymbol{z} \cdot \boldsymbol{z} \ a.e. \ in \ \Omega_2.$$

*However, for the sake of conciseness, we consider a scalar-valued coefficient.*

For a given source term $f \in L^2(\Omega)$, we consider the problem

$$\text{Find } u \in H_0^1(\Omega) \text{ such that } -\operatorname{div}(\varepsilon \nabla u) = f \in L^2(\Omega). \tag{1}$$

The equivalent variational formulation to (1) writes

$$\text{Find } u \in H_0^1(\Omega) \text{ such that } \int_\Omega \varepsilon \nabla u \cdot \nabla v \, \mathrm{d}\boldsymbol{x} = \int_\Omega f v \, \mathrm{d}\boldsymbol{x}, \quad \forall v \in H_0^1(\Omega). \tag{2}$$

Because of the change of sign of $\varepsilon$, the well-posedness of this problem does not fit into the classical theory of elliptic PDEs and it can be ill-posed. On the other hand, one can show that when $\kappa_\varepsilon^1$ or $\kappa_\varepsilon^2$ is large enough, Problem 2 is T-coercive (for instance see [6]), i.e. there exists an operator $T : H_0^1(\Omega) \to H_0^1(\Omega)$ such that $(u,v) \mapsto \int_\Omega \varepsilon \nabla u \cdot \nabla (T(v))$ is coercive, and then it is well-posed. For the case of polygonal interfaces, the construction of such operator T is based on the use of local isometric geometrical transformations (such as reflections, rotations, ...) near the interface, see [3].
The implementation of a general conforming finite element method to discretize (2) leads us to consider the problem

$$\text{Find } u_h \in V_h(\Omega) \text{ such that } \int_\Omega \varepsilon \nabla u_h \cdot \nabla v_h \, \mathrm{d}\boldsymbol{x} = \int_\Omega f v_h \, \mathrm{d}\boldsymbol{x}, \quad \forall v_h \in V_h(\Omega), \tag{3}$$

where $V_h(\Omega)$ is a well-chosen subspace $H_0^1(\Omega)$, and the parameter $h > 0$ is the so-called meshsize. Even in the case where (2) is T-coercive, one can not guaranty that Problem (3) is also T-coercive. Indeed, it may happen that for some $v_h \in V_h(\Omega)$, there holds $T(v_h) \notin V_h(\Omega)$. To overcome this difficulty, an interesting idea is to try to construct meshes such that the approximation spaces $V_h(\Omega)$ are stable by operators T for which Problem (2) is T-coercive. This type of meshes are called T-conform meshes. Such an approach has been investigated in [26, 12, 10]. It works quite well but presents two main drawbacks:

- The construction of well-suited meshes for curved interfaces, interfaces with corners or 3D interfaces is not an easy task [10, 3].

- Sometimes the operator T for which the problem is T-coercive is constructed by abstract tools and therefore is not explicit. In these situations, one cannot find adapted meshes.

On general meshes, three alternatives have already been proposed. The first one, was introduced in [6] and was based on the use of quasi-uniform meshes. In addition to this strong condition on the mesh, one of the limitation of this approach is that, for interfaces with general shapes, the convergence can not be assured in all the configurations in which Problem (2) is well-posed (because it is based on a particular (non-optimal) T-coercivity operator). The second one, presented in [14], consists in adding some dissipation to the problem (considering $\varepsilon + i\delta$ instead of $\varepsilon$ in (2) where $\delta$ depends on the meshsize). Unfortunately, this methods has a sub-optimal order of convergence even in the case where the solution and the interface is regular (see [14]). The third one is developed in [23] and is based on the use of mesh refinement techniques. Its essential limitation also lies the fact that, for interfaces with general shapes, the convergence can not be assured for all configurations in which Problem (2) is well-posed.

After that, in 2017, a new technique relying on the use of an optimal control reformulation has been introduced by Abdulle et al in [1]. Introducing

$$\mathrm{PH}^{1+s}(\Omega) := \{u \in H^1(\Omega)\,|\, u_{|\Omega_1} \in \mathrm{H}^{1+s}(\Omega_1) \text{ and } u_{|\Omega_2} \in \mathrm{H}^{1+s}(\Omega_2)\} \text{ for } s > 0,$$

their method is proved to be convergent for general meshes (that respect the interface) as soon as the exact solution to (1) belongs to the space $\mathrm{PH}^{1+s}(\Omega)$ for some $s > 1/2$. Unfortunately, this regularity condition is not always satisfied, especially when $\Sigma$ has corners in 2D or conical points in 3D. See the numerical illustration in Section 6.3 below.

In this work, we present a new strategy which relies on the use of a different optimal control reformulation and which converges without any restriction on the mesh (except the fact of being conforming to the interface), and without any restriction on the regularity of the exact solution. This method is inspired by the smooth extension method that was used (without proof of convergence) in [19] to approximate the solution to some classical scalar transmission problems.

The article is organized as follows. In Section 2, we start by giving a detailed description of the problem. Then, in Section 3, we explain how to derive an equivalent optimal control reformulation. Section 4 is dedicated to the study of some basic properties of the optimization problem and its regularization. The proposed numerical method and the proof of its convergence are given in Section 5. Our results are then illustrated by some numerical experiments in Section 6. Finally we give concluding remarks, including some possible extensions.

## 2  Main assumption on $\varepsilon$ and reformulation of the problem

Introduce the bounded operator $A_\varepsilon : \mathrm{H}_0^1(\Omega) \to (\mathrm{H}_0^1(\Omega))^*$ such that

$$\langle A_\varepsilon u, v \rangle = \int_\Omega \varepsilon \nabla u \cdot \nabla v \, \mathrm{d}\boldsymbol{x}, \qquad \forall u, v \in \mathrm{H}_0^1(\Omega).$$

Obviously $A_\varepsilon$ is an isomorphism if, and only if, Problem (1) is well-posed in the Hadamard sense. In this article, we shall work under the following

**Assumption 1.** *Assume that the coefficient $\varepsilon$ is such that $A_\varepsilon$ is an isomorphism.*

If $\varepsilon$ is piecewise constant by subdomain, the previous assumption is satisfied when the contrast $\kappa_\varepsilon := \varepsilon_2/\varepsilon_1$ does not belong to the so-called critical interval. The expression of this interval is in general not known analytically, except for particular geometries like symmetric domains, simple 2D interface with corners, simple 3D interfaces with circular conical tips (see [24, Chapter 2]). Under assumption 1, one is able to prove the accompanying shift theorem. We refer to [18, 7, 13, 12, 5].

**Theorem 2.1.** *Assume that $\Sigma$ is smooth (of class $\mathscr{C}^2$), polygonal (in 2D) or polyhedral (in 3D) and that Problem (1) is well-posed in the Hadamard sense. Then, there exists $\sigma_D(\varepsilon) \in (0,1]$ such that*

$$\forall f \in L^2(\Omega), \text{ the solution } u \text{ to Problem (1) is such that } \left| \begin{array}{ll} u \in \cup_{s \in [0, \sigma_D(\varepsilon))} PH^{1+s}(\Omega) & \text{if } \sigma_D(\varepsilon) < 1 \\ \\ u \in PH^2(\Omega) & \text{if } \sigma_D(\varepsilon) = 1 \end{array} \right.,$$

*with continuous dependence.*

The number $\sigma_D(\varepsilon)$ in the shift theorem is called the *(limit) regularity exponent*. For instance, when the interface is smooth and when it does not intersect with the boundary, then $\sigma_D(\varepsilon) = 1$ (cf. [18]).

**Remark 2.1.** *In Problem* (1), *we consider homogeneous Dirichlet boundary conditions. Let us mention that the results below extend quite straightforwardly to other situations, for example with Neumann or Robin-Fourier boundary conditions which can be homogeneous or not, as long as the associated operator is an isomorphism.*

To introduce the method, we start by writing an equivalent version of (1) in which the unknown $u \in H_0^1(\Omega)$ is split into two unknowns defined in $\Omega_1$ and $\Omega_2$ : $(u_1, u_2) := (u_{|\Omega_1}, u_{|\Omega_2})$. To do so, we observe that since $f \in L^2(\Omega)$, the solution $u$ to (1) is such that the vector field $\varepsilon \nabla u$ belongs to the space $H(\text{div}, \Omega) = \{\boldsymbol{u} \in (L^2(\Omega))^d \text{ such that } \text{div}(\boldsymbol{u}) \in L^2(\Omega)\}$. Consequently, the pair of functions $(u_1, u_2)$ satisfies the problem

$$\text{Find } (u_1, u_2) \in V_1(\Omega_1) \times V_2(\Omega_2) \text{ such that } \begin{vmatrix} -\text{div}(\varepsilon_1 \nabla u_1) = f_1 =: f_{|\Omega_1} \\ -\text{div}(\varepsilon_2 \nabla u_2) = f_2 =: f_{|\Omega_2} \\ \varepsilon_1 \partial_n u_1 = \varepsilon_2 \partial_n u_2 \text{ and } u_1 = u_2 \text{ on } \Sigma \end{vmatrix} \qquad (4)$$

in which $\boldsymbol{n}$ stands for the unit normal vector to $\Sigma$ oriented to the exterior of $\Omega_2$ (see Figure 1) and

$$V_1(\Omega_1) := \{u \in H^1(\Omega_1), u = 0 \text{ on } \partial\Omega_1 \backslash \Sigma\}, \qquad V_2(\Omega_2) := \{u \in H^1(\Omega_2), u = 0 \text{ on } \partial\Omega_2 \backslash \Sigma\}.$$

Since $\text{meas}_{\partial\Omega}(\partial\Omega_2 \setminus \Sigma) > 0$, all elements of $V_2(\Omega_2)$ fulfill a homogeneous boundary condition on a part of the boundary $\partial\Omega_2$. On the other hand, one can check that if $(u_1, u_2)$ is a solution to (4), then the function $u$ defined by $u_{|\Omega_j} = u_i$ for $j = 1, 2$ solves (1). The equations satisfied by $u_1$ and $u_2$ are elliptic but they are coupled by the transmission conditions on $\Sigma$. As a consequence, we cannot solve them independently. The purpose of the next paragraph is to explain how to proceed to write an alternative formulation (an optimization-based one), which can be solved via an iterative procedure such that at each step one has to solve a set of elliptic problems.

## 3 The smooth extension method and optimal control reformulation of the problem

The smooth extension method was proposed in [21] and can be considered as a special case of the fictitious domain methods (see [2]). It has been adapted to study the classical scalar transmission problem, i.e. with constant sign coefficients, in [19]. In this section, we explain how to apply it to our problem.

### 3.1 Presentation of the smooth extension method

The idea behind the smooth extension method is the following: instead of looking for $(u_1, u_2) \in V_1(\Omega) \times V_2(\Omega_2)$ solution to (4), we search for a pair of functions $(\tilde{u}, u_2) \in H_0^1(\Omega) \times V_2(\Omega_2)$ such that $(\tilde{u}_{|\Omega_1}, u_2)$ is a solution to (4).[1] The function $\tilde{u}$ is then a particular continuous extension of $u_1$ to the whole domain $\Omega$. The difficulty is to find a "good" way to define the function $\tilde{u}$ so that it can be approximated by the classical FEM. The function $u_2$ can then be approximated by solving the elliptic problem satisfied by $u_2$ in $\Omega_2$ completed by $\tilde{u}_{|\Sigma}$ (resp. $\varepsilon_1 \partial_n \tilde{u}_{|\Sigma}$) as a Dirichlet (resp. Neumann) boundary condition on $\Sigma$. Note that at first sight the construction of such $\tilde{u}$ is not straightforward. This will be achieved thanks to an optimal control reformulation of (4). This is the main goal of the next paragraph in which we also reformulate the idea presented above in a more rigorous way.

---

[1]In the text below, we choose an extension from $\Omega_1$ to $\Omega^\star = \Omega_2$. Obviously, one could choose an extension from $\Omega_2$ to $\Omega^\star = \Omega_1$ so that $(u_1, \tilde{u}_{|\Omega_2})$ is a solution to (4). In this case, the condition $\text{meas}_{\partial\Omega}(\partial\Omega_1 \setminus \Sigma) > 0$ must hold.

## 3.2 An optimal control reformulation of the problem

Before getting into details, let us first introduce $\tilde{\varepsilon}_1 \in \mathrm{L}^\infty(\Omega)$ such that $\tilde{\varepsilon}_1 \geq \tilde{\varepsilon}^+ > 0$ a. e. in $\Omega$ and $\tilde{\varepsilon}_1 = \varepsilon_1$ in $\Omega_1$. Then, let $E : \mathrm{V}_1(\Omega_1) \to \mathrm{H}_0^1(\Omega)$ be an arbitrary continuous extension operator. By making use of (4), one can show easily that

$$\left| \begin{array}{ll} \int_\Omega \tilde{\varepsilon}_1 \nabla E(u_1) \cdot \nabla v \, \mathrm{d}\boldsymbol{x} = \int_{\Omega_1} f_1 v \, \mathrm{d}\boldsymbol{x} + \int_{\Omega_2} \tilde{\varepsilon}_1 \nabla E(u_1) \cdot \nabla v \, \mathrm{d}\boldsymbol{x} - \langle \varepsilon_1 \partial_n u_1, v \rangle_\Sigma & \forall v \in \mathrm{H}_0^1(\Omega), \\ \int_{\Omega_2} \varepsilon_2 \nabla u_2 \cdot \nabla v_2 \, \mathrm{d}\boldsymbol{x} = \int_{\Omega_2} f_2 \, v_2 \, \mathrm{d}\boldsymbol{x} + \langle \varepsilon_1 \partial_n u_1, v_2 \rangle_\Sigma & \forall v_2 \in \mathrm{V}_2(\Omega_2). \end{array} \right.$$

Now, given that the linear form $v_2 \mapsto \int_{\Omega_2} \tilde{\varepsilon}_1 \nabla E(u_1) \cdot \nabla v_2 \, \mathrm{d}\boldsymbol{x} - \langle \varepsilon_1 \partial_n u_1, v_2 \rangle_\Sigma$ is continuous on $\mathrm{V}_2(\Omega_2)$ one can define, thanks to the Riesz representation theorem, for each $E(u_1)$ a unique $w_{E(u_1)} \in \mathrm{V}_2(\Omega_2)$ such that

$$\int_{\Omega_2} \tilde{\varepsilon}_1 \nabla E(u_1) \cdot \nabla v_2 \, \mathrm{d}\boldsymbol{x} - \langle \varepsilon_1 \partial_n u_1, v_2 \rangle_\Sigma = \int_{\Omega_2} \tilde{\varepsilon}_1 \nabla w_{E(u_1)} \cdot \nabla v_2 \, \mathrm{d}\boldsymbol{x} \quad \forall v_2 \in \mathrm{V}_2(\Omega_2). \tag{5}$$

Above we have used the fact that $(u,v) \mapsto (\tilde{\varepsilon}_1 \nabla u, \nabla v)_{\mathrm{L}^2(\Omega)^d}$ is an inner product on $\mathrm{V}_2(\Omega_2)$. As a consequence, we have

$$\left| \begin{array}{ll} \int_\Omega \tilde{\varepsilon}_1 \nabla E(u_1) \cdot \nabla v \, \mathrm{d}\boldsymbol{x} = \int_{\Omega_1} f_1 v \, \mathrm{d}\boldsymbol{x} + \int_{\Omega_2} \tilde{\varepsilon}_1 \nabla w_{E(u_1)} \cdot \nabla v \, \mathrm{d}\boldsymbol{x} & \forall v \in \mathrm{H}_0^1(\Omega), \\ \int_{\Omega_2} \varepsilon_2 \nabla u_2 \cdot \nabla v_2 \, \mathrm{d}\boldsymbol{x} = \int_{\Omega_2} f_2 \, v_2 \, \mathrm{d}\boldsymbol{x} + \int_{\Omega_2} \tilde{\varepsilon}_1 \nabla (E(u_1) - w_{E(u_1)}) \cdot \nabla v_2 \, \mathrm{d}\boldsymbol{x} & \forall v_2 \in \mathrm{V}_2(\Omega_2). \end{array} \right.$$

Since the coefficients $\tilde{\varepsilon}_1$ and $\varepsilon_2$ have fixed signs, the forms

$$(u,v) \mapsto \int_\Omega \tilde{\varepsilon}_1 \nabla u \cdot \nabla v \, \mathrm{d}\boldsymbol{x} \text{ and } (u_2, v_2) \mapsto -\int_{\Omega_2} \varepsilon_2 \nabla u_2 \cdot \nabla v_2 \, \mathrm{d}\boldsymbol{x},$$

are coercive, respectively on $\mathrm{H}_0^1(\Omega)$ and on $\mathrm{V}_2(\Omega_2)$. With this in mind, we define for all $w \in \mathrm{V}_2(\Omega_2)$, the couple of functions $(u^w, u_2^w) \in \mathrm{H}_0^1(\Omega) \times \mathrm{V}_2(\Omega_2)$ that are solution to the well-posed system of equations:

$$\left| \begin{array}{ll} \int_\Omega \tilde{\varepsilon}_1 \nabla u^w \cdot \nabla v \, \mathrm{d}\boldsymbol{x} = \int_{\Omega_1} f_1 v \, \mathrm{d}\boldsymbol{x} + \int_{\Omega_2} \tilde{\varepsilon}_1 \nabla w \cdot \nabla v \, \mathrm{d}\boldsymbol{x} & \forall v \in \mathrm{H}_0^1(\Omega), \\ \int_{\Omega_2} \varepsilon_2 \nabla u_2^w \cdot \nabla v_2 \, \mathrm{d}\boldsymbol{x} = \int_{\Omega_2} f_2 \, v_2 \, \mathrm{d}\boldsymbol{x} + \int_{\Omega_2} \tilde{\varepsilon}_1 \nabla (u^w - w) \cdot \nabla v_2 \, \mathrm{d}\boldsymbol{x} & \forall v_2 \in \mathrm{V}_2(\Omega_2). \end{array} \right. \tag{6}$$

Well-posedness is achieved by solving the elliptic problem in $u^w$ first, and then the elliptic problem in $u_2^w$.

**Remark 3.1.** *Let $w \in \mathrm{V}_2(\Omega)$ be constructed by means of a particular continuous extension via (5) of $u_1$, the first part of the solution to (4): that is, with $w = w_{E(u_1)}$. Then it follows that $(u^w_{|\Omega_1}, u_2^w)$ is the solution to (4). Indeed, one finds first that $u^w = E(u_1)$, and then that $u_2^w = u_2$.*

On the other hand, we have the following result

**Proposition 3.1.** *For all $w \in \mathrm{V}_2(\Omega_2)$, the functions $u_1^w := u^w_{|\Omega_1}$ and $u_2^w$ are such that*

$$\left| \begin{array}{ll} -\mathrm{div}(\varepsilon_1 \nabla u_1^w) = f_1 & \text{in } \Omega_1, \\ -\mathrm{div}(\varepsilon_2 \nabla u_2^w) = f_2 & \text{in } \Omega_2 \\ \varepsilon_1 \partial_n u_1^w = \varepsilon_2 \, \partial_n u_2^w & \text{on } \Sigma. \end{array} \right.$$

*Proof.* Take $\varphi_1 \in \mathscr{C}_0^\infty(\Omega_1)$ and extend it by 0 to the whole $\Omega$ to obtain the function $\varphi \in \mathscr{C}_0^\infty(\Omega)$. Take $v = \varphi$ in the problem satisfied by $u^w$. One finds that $-\text{div}(\varepsilon_1 \nabla u_1^w) = f_1$ in $\Omega_1$. Next, take some $\varphi_2 \in \mathscr{C}_0^\infty(\Omega_2)$, extend it by 0 in $\Omega_1$ and denote by $\varphi$ the new function. By taking $v = \varphi$ in the problem satisfied by $u^w$ and $v_2 = \varphi_2$ in the problem satisfied by $u_2^w$ one finds that

$$\int_{\Omega_2} \tilde{\varepsilon}_1 \nabla u^w \cdot \nabla \varphi_2 \, \mathrm{d}\boldsymbol{x} = \int_{\Omega_2} \tilde{\varepsilon}_1 \nabla w \cdot \nabla \varphi_2 \, \mathrm{d}\boldsymbol{x},$$
$$\int_{\Omega_2} \varepsilon_2 \nabla u_2^w \cdot \nabla \varphi_2 \, \mathrm{d}\boldsymbol{x} = \int_{\Omega_2} f_2 \varphi_2 \, \mathrm{d}\boldsymbol{x} + \int_{\Omega_2} \tilde{\varepsilon}_1 \nabla (u^w - w) \cdot \nabla \varphi_2 \, \mathrm{d}\boldsymbol{x}.$$

By considering the sum of the two formulations, we conclude that $-\text{div}(\varepsilon_2 \nabla u_2^w) = f_2$ in $\Omega_2$. To end the proof, it remains to show that $\varepsilon_1 \partial_n u^w = \varepsilon_2 \partial_n u_2^w$. For this, let $v \in \mathrm{H}_0^1(\Omega)$ and define $v_2 = v_{|\Omega_2} \in \mathrm{V}_2(\Omega_2)$. By taking $v$ and $v_2$ as test functions in (6), considering the sum of the two equations, integrating by parts in both formulations and then, using the equations satisfied by $u_1^w$ and $u_2^w$, we infer that

$$-\langle \varepsilon_1 \partial_n u_1^w, v \rangle = -\langle \varepsilon_2 \partial_n u_2^w, v \rangle, \qquad v \in \mathrm{H}_0^1(\Omega).$$

According to the surjectivity of the trace mapping on $\Sigma$, this gives $\varepsilon_1 \partial_n u_1^w = \varepsilon_2 \partial_n u_2^w$ on $\Sigma$ and ends the proof. $\square$

Thus the introduction of an auxiliary "control" function $w \in \mathrm{V}_2(\Omega_2)$ allows us to construct pseudo-solutions to the equation (4) for which the condition on the normal derivatives is automatically satisfied. However we do not have in general continuity across interface. Taking this into account, we get the

**Lemma 3.1.** *If there exists $w^* \in \mathrm{V}_2(\Omega_2)$ such that the solution to (6) satisfies $u^{w^*} = u_2^{w^*}$ on $\Sigma$, then $(u^{w^*}_{|\Omega_1}, u_2^{w^*})$ solves (4).*

Thanks to what we have explained in Remark 3.1, we know that to every continuous extension of $u_1$ to $\Omega$, one can define $w^* \in \mathrm{V}_2(\Omega_2)$ for which $u^{w^*} = u_2^{w^*}$ on $\Sigma$. This leads us to the following result.

**Lemma 3.2.** *Let $u_1$ be the first part of the solution to (4). Then, the set of $w^* \in \mathrm{V}_2(\Omega_2)$ such that the solution to (6) satisfies $u^{w^*} = u_2^{w^*}$ on $\Sigma$ is isomorphic to the set of all possible continuous extensions of $u_1$ to $\Omega$. Furthermore, $w^*$ and $u^{w^*}$ are linked by relation*

$$\int_{\Omega_2} \tilde{\varepsilon}_1 \nabla u^{w^*} \cdot \nabla v_2 \, \mathrm{d}\boldsymbol{x} - \langle \varepsilon_1 \partial_n u_1, v_2 \rangle_\Sigma = \int_{\Omega_2} \tilde{\varepsilon}_1 \nabla w^* \cdot \nabla v_2 \, \mathrm{d}\boldsymbol{x} \quad \forall v_2 \in \mathrm{V}_2(\Omega_2). \tag{7}$$

Now, we have all the tools to introduce the optimal control reformulation of the problem (4). As a matter of fact, in order to find a function $w^* \in \mathrm{V}_2(\Omega_2)$ for which $u^{w^*} = u_2^{w^*}$ on $\Sigma$, it is enough to solve the following optimal control problem:

$$\text{Find } w^* = \underset{w \in \mathrm{V}_2(\Omega_2)}{\operatorname{argmin}} J(w) \quad \text{with} \quad J(w) = \frac{1}{2} \int_\Sigma |u^w - u_2^w|^2 \, \mathrm{d}\sigma, \tag{8}$$

where $(u^w, u_2^w) \in \mathrm{H}_0^1(\Omega) \times \mathrm{V}_2(\Omega_2)$ is the solution to (6). Note that in (8), the functional $J$ plays the role of the cost functional, while (6) plays the role of the state equation. Obviously, thanks to Lemma 3.2, the problem (8) has an infinite number of solutions. As a result, one may need to use a regularization technique in order to be able to construct a convergent discretization of the problem: this will be the subject of §4.3 where we will study the classical Tikhonov regularization method applied to Problem (8).

# 4 Basic properties of the optimization problem and its regularization

In this section, we present in §4.1 some useful properties of the cost functional $J$ and of the set of its minimizers in §4.2 . After that in §4.3, we study the Tikhonov regularization of the problem. Furthermore, we explain, in §4.4, how to use the the adjoint approach in order to find an explicit expression of the gradient of $J$.

## 4.1 Properties of the cost functional

Since we have used the $\mathrm{L}^2(\Sigma)$ norm instead of the $\mathrm{H}^{1/2}(\Sigma)$ norm in the definition of $J$, one has the following results.

**Proposition 4.1.** *The cost functional $J$ satisfies the following properties:*

1. *Let $(w_n)_n$ be a sequence of elements of $\mathrm{V}_2(\Omega_2)$ that converges weakly to $w_0 \in \mathrm{V}_2(\Omega_2)$. Then, $(J(w_n))_n$ converges to $J(w_0)$.*

2. *The functional $J$ is continuous and convex on $\mathrm{V}_2(\Omega_2)$.*

*Proof.* 1. For all $n \in \mathbb{N}$, denote by $(u^n, u_2^n) \in \mathrm{H}_0^1(\Omega) \times \mathrm{V}_2(\Omega)$ the solution to (6) with $w = w_n$. From the ellipticity of the problems involved in (6), it follows that $(u^n)_n$ (resp. $(u_2^n)_n$) converges weakly in $\mathrm{H}_0^1(\Omega)$ (resp. $\mathrm{V}_2(\Omega_2)$) to some $u \in \mathrm{H}^1(\Omega)$ (resp. $u_2 \in \mathrm{V}_2(\Omega_2)$) such that $(u, u_2)$ is the solution to (6) with $w = w_0$. The continuity of the trace operator from $\mathrm{H}^1(\Omega)$ to $\mathrm{H}^{1/2}(\Sigma)$ ensures that $(u_{|\Sigma}^n - u_{2|\Sigma}^n)_n$ converges weakly to $u_{|\Sigma} - u_{2|\Sigma}$ in $\mathrm{H}^{1/2}(\Sigma)$. Given that the embedding of $\mathrm{H}^{1/2}(\Sigma)$ into $\mathrm{L}^2(\Sigma)$ is compact, it actually converges strongly to $u_{|\Sigma} - u_{2|\Sigma}$ in $\mathrm{L}^2(\Sigma)$. Thus $(J(w_n))$ converges to $J(w_0)$. The result is proved.

2. While the continuity is a direct consequence of the first item, the convexity follows from the fact that $J : \mathrm{V}_2(\Omega_2) \to \mathbb{R}$ is the composition of the affine map $j_1 : \mathrm{V}_2(\Omega_2) \to \mathrm{L}^2(\Sigma)$ and of the convex map $j_2 : \mathrm{L}^2(\Sigma) \to \mathbb{R}$ such that for all $w \in \mathrm{V}_2(\Omega_2)$, $g \in \mathrm{L}^2(\Sigma)$ we have

$$\left| \begin{aligned} &j_1(w) = (u^w - u_2^w)_{|\Sigma} \text{ where } (u^w, u_2^w) \in \mathrm{H}_0^1(\Omega) \times \mathrm{V}_2(\Omega_2) \text{ is the solution to (6)}, \\ &j_2(g) = \frac{1}{2} \int_\Sigma |g|^2 \, \mathrm{d}\sigma. \end{aligned} \right. \qquad (9)$$

□

## 4.2 The set of minimizers of the functional $J$

Thanks to Lemma 3.2, we know that $J$ has an infinite number of minimizers. This (non-empty) set will be denoted by $M_J$. Without any difficulty, one can see that $M_J$ coincides with the set of zeros of the functional $J$. As a result, since $J$ is continuous, convex and positive, the set $M_J$ is closed and convex in $\mathrm{V}_2(\Omega_2)$. This allows us to say that the following minimization problem:

$$\min_{w \in M_J} \int_{\Omega_2} \tilde{\varepsilon}_1 |\nabla w|^2 \, \mathrm{d}\boldsymbol{x}$$

has a unique solution, as a consequence of the strict convexity of $v_2 \mapsto \int_{\Omega_2} \tilde{\varepsilon}_1 |\nabla v_2|^2 \, \mathrm{d}\boldsymbol{x}$ in $\mathrm{V}_2(\Omega_2)$, and of the fact that $M_J$ is a closed, convex subset of $\mathrm{V}_2(\Omega_2)$. In the following, we shall denote by $w_J^*$ the smallest minimizer of the functional $J$:

$$w_J^* = \underset{w \in M_J}{\operatorname{argmin}} \int_{\Omega_2} \tilde{\varepsilon}_1 |\nabla w|^2 \, \mathrm{d}\boldsymbol{x}. \qquad (10)$$

The goal of the rest of this paragraph is to find a characterization of $E_{w_J^*}(u_1)$, the continuous extension of $u_1$ that is associated with $w_J^*$. Note that the link between $E_{w_J^*}(u_1)$ and $w_J^*$ is given by the following (see relation (7)):

$$\int_{\Omega_2} \tilde{\varepsilon}_1 \nabla E_{w_J^*}(u_1) \cdot \nabla v_2 \, \mathrm{d}\boldsymbol{x} - \langle \varepsilon_1 \partial_n u_1, v_2 \rangle_\Sigma = \int_{\Omega_2} \tilde{\varepsilon}_1 \nabla w_J^* \cdot \nabla v_2 \, \mathrm{d}\boldsymbol{x} \quad \forall v_2 \in \mathrm{V}_2(\Omega_2). \quad (11)$$

To proceed, we define $E_H(u_1) \in \mathrm{H}_0^1(\Omega)$ the continuous extension of $u_1$ that satisfies

$$\mathrm{div}(\tilde{\varepsilon}_1 \nabla E_H(u_1)) = 0 \text{ in } \Omega_2. \quad (12)$$

In particular, we have

$$\int_{\Omega_2} \tilde{\varepsilon}_1 \nabla E_H(u_1) \cdot \nabla v_2 \, \mathrm{d}\boldsymbol{x} = 0 \quad \forall v_2 \in \mathrm{H}_0^1(\Omega_2).$$

Denote by $w_H \in M_J$ the minimizer associated with $E_H(u_1)$. Thanks to (7), we know that

$$\int_{\Omega_2} \tilde{\varepsilon}_1 \nabla E_H(u_1) \cdot \nabla v_2 \, \mathrm{d}\boldsymbol{x} - \langle \varepsilon_1 \partial_n u_1, v_2 \rangle_\Sigma = \int_{\Omega_2} \tilde{\varepsilon}_1 \nabla w_H \cdot \nabla v_2 \, \mathrm{d}\boldsymbol{x} \quad \forall v_2 \in \mathrm{V}_2(\Omega_2). \quad (13)$$

We infer that

$$\int_{\Omega_2} \tilde{\varepsilon}_1 \nabla w_H \cdot \nabla v_2 \, \mathrm{d}\boldsymbol{x} = 0 \quad \forall v_2 \in \mathrm{H}_0^1(\Omega_2). \quad (14)$$

By taking the difference between (11) and (13), taking $v_2 = w_H$, using the fact that $E_H(u_1) - E_{w_J^*}(u_1) \in \mathrm{H}_0^1(\Omega_2)$ and owing to (14), we infer that

$$\int_{\Omega_2} \tilde{\varepsilon}_1 \nabla(w_H - w_J^*) \cdot \nabla w_H \, \mathrm{d}\boldsymbol{x} = 0,$$

so

$$\int_{\Omega_2} \tilde{\varepsilon}_1 |\nabla w_J^*|^2 \, \mathrm{d}\boldsymbol{x} = \int_{\Omega_2} \tilde{\varepsilon}_1 |\nabla w_J^* - \nabla w_H|^2 \, \mathrm{d}\boldsymbol{x} + \int_{\Omega_2} \tilde{\varepsilon}_1 |\nabla w_H|^2 \, \mathrm{d}\boldsymbol{x} \geq \int_{\Omega_2} \tilde{\varepsilon}_1 |\nabla w_H|^2 \, \mathrm{d}\boldsymbol{x}.$$

Hence, from the definition of $w_J^*$, we then obtain the following

**Proposition 4.2.** *The functions $w_H$ and $w_J^*$ coincide.*

**Remark 4.1.** *It is worth noting that, thanks to (13) and using the definition of $E_H(u_1)$, the function $w_H$ satisfies the problem:*

$$\mathrm{div}(\tilde{\varepsilon} \nabla w_H) = 0 \text{ in } \Omega_2 \text{ and } \tilde{\varepsilon}_1 \partial_n w_{H|\Sigma} = \tilde{\varepsilon}_1 \partial_n E_H(u_1)_{|\Sigma} - \varepsilon_1 \partial_n u_{1|\Sigma}. \quad (15)$$

*Recall that $\boldsymbol{n}$ is the unit normal vector to $\Sigma$ oriented to the exterior of $\Omega_2$.*

### 4.3 Tikhonov regularization of the problem

Tikhonov regularization, which was originally introduced in [25], is a classical method to regularize a convex optimization problem. Classically, this method is used in the context of regularization of ill-posed inverse problems (see [20] and the references therein). In this paragraph, we study the convergence of such regularization when it is applied to our problem. For $\delta > 0$, we introduce the functional $J^\delta : \mathrm{V}_2(\Omega_2) \to \mathbb{R}$ defined by

$$J^\delta(w) = J(w) + \delta \int_{\Omega_2} \tilde{\varepsilon}_1 |\nabla w|^2 \, \mathrm{d}\boldsymbol{x} \quad \forall w \in \mathrm{V}_2(\Omega_2).$$

To simplify notation, we will denote by $\|\cdot\|_{\tilde{\varepsilon}_1} : V_2(\Omega_2) \to \mathbb{R}_+$ the norm that is defined as follows:

$$\|w\|_{\tilde{\varepsilon}_1} := \left(\int_{\Omega_2} \tilde{\varepsilon}_1 |\nabla w|^2 \, d\boldsymbol{x}\right)^{1/2}, \quad \forall w \in V_2(\Omega_2).$$

Endowed with the associated scalar product $(\cdot,\cdot)_{\tilde{\varepsilon}_1}$, the space $V_2(\Omega_2)$ is a Hilbert space. Since $J$ is convex and $\delta > 0$, the functional $J^\delta$ is strictly convex and coercive. Therefore the minimization problem

$$\min_{w \in V_2(\Omega_2)} J^\delta(w)$$

has a unique solution that we denote by $w_\delta^*$. Our goal is to study the behaviour of $w_\delta^*$ as $\delta$ tends to zero. One expects $(w_\delta^*)_\delta$ to converge to one of the solutions (8). If this is the case and because the problem (8) has an infinite number of solutions, it will be interesting to characterize the particular solution towards which $(w_\delta^*)_\delta$ converges. Our findings are summarized in the following

**Proposition 4.3.** *When $\delta \to 0$, the sequence $(w_\delta^*)_\delta$ converges towards $w_J^*$, the smallest minimizer of $J$.*

The proof of the previous result is quite classical. However, for the convenience of the reader, we will detail it in Appendix A.
In conclusion, we can say that the Tikhonov regularization allows us to obtain a stabilized version of the optimization problem (8). This will be used in order to introduce a stabilization of the discretization of the problem (8), but in that case the stabilization parameter $\delta$ will be chosen as a function of the discretization parameter. This will be detailed in §5.3. Note that the same idea was employed in [1].

### 4.4 Gradient of the functional $J$

As indicated in the introduction, the main objective of this work is to propose a new numerical method for approximating the solution to (1). This method will be based on the numerical approximation of the solution to the optimization problem (8). In this section, we will explain how to obtain an explicit expression of $J'(w)$ the gradient of $J$ at some $w \in V_2(\Omega)$. We recall that the functional $J$ is differentiable, because it can be written as a composition of the two differentiable maps $j_1$ and $j_2$, cf. (9). And, since the functional $J$ is scalar valued, its differential at $w \in V_2(\Omega_2)$ can be represented by its gradient $J'(w) \in V_2(\Omega_2)$:

$$\text{For all } h \in V_2(\Omega_2), \quad \int_{\Omega_2} \tilde{\varepsilon}_1 \nabla J'(w) \cdot \nabla h \, d\boldsymbol{x} = \lim_{t \to 0} \frac{J(w + th) - J(w)}{t}.$$

To find an explicit expression of $J'(w)$, we use the adjoint approach [11]. Details about the application of this approach to our problem are given in Appendix B (see also [19]). Here, we present final result. To do so, we start by introducing the so-called adjoint equations. For all $w \in V_2(\Omega_2)$, recalling that $(u^w, u_2^w) \in H_0^1(\Omega) \times V_2(\Omega_2)$ is the solution to (6), we introduce $(g^w, g_2^w) \in H_0^1(\Omega) \times V_2(\Omega_2)$ such that

$$\left|\begin{array}{ll} \int_{\Omega} \tilde{\varepsilon}_1 \nabla g^w \cdot \nabla v \, d\boldsymbol{x} = \int_{\Omega_2} \tilde{\varepsilon}_1 \nabla g_2^w \cdot \nabla v \, d\boldsymbol{x} - \int_{\Sigma} (u^w - u_2^w) v \, d\sigma & \forall v \in H_0^1(\Omega) \\ \int_{\Omega_2} \varepsilon_2 \nabla g_2^w \cdot \nabla v_2 \, d\boldsymbol{x} = \int_{\Sigma} (u^w - u_2^w) v_2 \, d\sigma & \forall v_2 \in V_2(\Omega_2). \end{array}\right. \quad (16)$$

As observed before, the functions $g^w$, $g_2^w$ are well-defined. One can prove the

**Lemma 4.1.** *For all $w \in V_2(\Omega_2)$, there holds $J'(w) = g_2^w - g^w_{|\Omega_2}$, where $(g^w, g_2^w)$ solve (16).*

We have the following optimality result

**Corollary 4.1.** *We have the equivalence*

$$[\, w^* \in V_2(\Omega_2) \text{ is such that } J'(w^*) = 0 \,] \quad \Longleftrightarrow \quad w^* \in M_J.$$

*Proof.* Let us start with the proof of the direct implication. Suppose that there exists some $w^* \in V_2(\Omega_2)$ such that $g^{w^*}_{|\Omega_2} = g_2^{w^*}$. By taking the sum of the variational formulations of (16), we deduce that

$$\int_\Omega \varepsilon \nabla g^{w^*} \cdot \nabla v \, d\boldsymbol{x} = 0 \qquad \forall v \in H^1_0(\Omega).$$

This means $A_\varepsilon(g^{w^*}) = 0$ and then, thanks to Assumption 1, $g^{w^*} = 0$. This implies that $g_2^{w^*} = 0$ and then by using the second equation of (16), that $u^{w^*} = u_2^{w^*}$ on $\Sigma$. This shows that $w^*$ is a minimizer of $J$. The reverse implication is a consequence of the fact that if $w^* \in M_J$ we have $J(w^*) = 0$ and then $u^{w^*} = u_2^{w^*}$ on $\Sigma$. This implies that $g_2^{w^*} = 0$ and that $g^{w^*} = 0$. □

We end this paragraph with the following result that can be useful to prove the convergence of the classical gradient descent algorithm.

**Corollary 4.2.** *The functional $J' : V_2(\Omega_2) \to V_2(\Omega_2)$ is Lipschitz continuous.*

*Proof.* Starting from (6), we deduce that $w \mapsto u^w$, $w \mapsto u_2^w$ are Lipschitz continuous. Inserting this into (16), we obtain the result. □

## 5 Numerical discretization of the problem

In this part, we are concerned with the numerical approximation of (8) by means of the Finite Element Method. To do so, we start by presenting some details and notations about the family of meshes that will be used. To simplify the presentation, the domain $\Omega$ and the subdomains $(\Omega_i)_{i=1,2}$ are supposed to have polygonal (when $d = 2$) or polyhedral (resp. $d = 3$) boundaries.

### 5.1 Meshes and discrete spaces

Let $(\mathcal{T}_h)_h$ be a regular family of meshes of $\overline{\Omega}$ (see [15]), composed of (closed) simplices. The subscript $_h$ stands for the meshsize.

**Assumption 2.** *We suppose that for all $h$, every simplex of $\mathcal{T}_h$ belongs either to $\overline{\Omega_1}$ or to $\overline{\Omega_2}$.*

According to Assumption 2, for $i = 1, 2$, one can consider the family of meshes $(\mathcal{T}_h^i)_h$ made of those simplices that belong to $\overline{\Omega_i}$.
For all $k \in \mathbb{N}^*$, we set

$$V_h^k(\Omega) := \{v_h \in H^1_0(\Omega) \,|\, v_{h|T} \in P^k(T), \, \forall T \in \mathcal{T}_h\}.$$

Here $P^k(T)$ stands for the space of polynomials (of $d$ variables) defined on $T$ of degree at most equal to $k$. In the same way, we define for $i = 1, 2$,

$$V_h^k(\Omega_i) := \{v_{i,h} \in H^1(\Omega_i) \,|\, v_{i,h|T} \in P^k(T), \forall T \in \mathcal{T}_h^i \text{ and } v_{i,h} = 0 \text{ on } \partial \Omega_i \backslash \Sigma\}.$$

**Remark 5.1.** *Note that under Assumption 2 for all $h > 0$, for $i = 1, 2$ the space $V_h^k(\Omega_i)$ coincides with $\{u_{|\Omega_i} \,|\, u \in V_h^k(\Omega)\}$.*

Finally, we recall the basic approximability properties

$$\left| \begin{array}{l} \forall v \in H^1_0(\Omega), \quad \lim_{h \to 0} \left( \inf_{v_h \in V_h^k(\Omega)} \|v - v_h\|_{H^1_0(\Omega)} \right) = 0, \\ \forall v_2 \in V_2(\Omega_2), \quad \lim_{h \to 0} \left( \inf_{v_{2,h} \in V_h^k(\Omega_2)} \|v_2 - v_{2,h}\|_{\tilde{\varepsilon}_1} \right) = 0. \end{array} \right. \qquad (17)$$

## 5.2 Discretization strategy

For $h > 0$ and $w \in V_2(\Omega)$, define the functions $u_h^w \in V_h^k(\Omega)$ and $u_{2,h}^w \in V_h^k(\Omega_2)$ as the solutions to the following well-posed discrete problems:

$$\left| \begin{array}{l} \int_\Omega \tilde{\varepsilon}_1 \nabla u_h^w \cdot \nabla v_h \, \mathrm{d}\boldsymbol{x} = \int_{\Omega_1} f v_h \, \mathrm{d}\boldsymbol{x} + \int_{\Omega_2} \tilde{\varepsilon}_1 \nabla w \cdot \nabla v_h \, \mathrm{d}\boldsymbol{x}, \; \forall v_h \in V_h^k(\Omega) \\ \int_{\Omega_2} \varepsilon_2 \nabla u_{2,h}^w \cdot \nabla v_{2,h} \, \mathrm{d}\boldsymbol{x} = \int_{\Omega_2} f_2 v_{2,h} \, \mathrm{d}\boldsymbol{x} + \int_{\Omega_2} \tilde{\varepsilon}_1 \nabla (u_h^w - w) \cdot \nabla v_{2,h} \, \mathrm{d}\boldsymbol{x}, \; \forall v_{2,h} \in V_h^k(\Omega_2). \end{array} \right. \tag{18}$$

Then introduce the projection operator $\pi_h^k : V_2(\Omega_2) \to V_h^k(\Omega_2)$ such that for all $w \in V_2(\Omega_2)$, $\pi_h^k w$ is defined as the unique element of $V_h^k(\Omega_2)$ that satisfies the problem

$$\int_{\Omega_2} \tilde{\varepsilon}_1 \nabla (\pi_h^k w) \cdot \nabla v_{2,h} \, \mathrm{d}\boldsymbol{x} = \int_{\Omega_2} \tilde{\varepsilon}_1 \nabla w \cdot \nabla v_{2,h} \, \mathrm{d}\boldsymbol{x} \quad \forall v_{2,h} \in V_h^k(\Omega_2).$$

Obviously, one has the estimate

$$\|\pi_h^k w\|_{\tilde{\varepsilon}_1} \leq \|w\|_{\tilde{\varepsilon}_1}. \tag{19}$$

From the definition of $\pi_h^k w$, one can easily see that for all $w \in V_2(\Omega_2)$ we have the identities

$$u_h^{\pi_h^k w} = u_h^w \qquad \text{and} \qquad u_{2,h}^{\pi_h^k w} = u_{2,h}^w. \tag{20}$$

Now, let us turn our attention to the discretization of the optimization problem (8). The natural way to do that is to replace it by the problem

$$\inf_{w_h \in V_h^k(\Omega_2)} J_{0,h}(w_h) := \frac{1}{2} \int_\Sigma |u_h^{w_h} - u_{2,h}^{w_h}|^2 \, \mathrm{d}\sigma. \tag{21}$$

One can proceed as in the proof of proposition 4.1 to show that the cost functional $J_{0,h} : V_h^k \to \mathbb{R}$ (defined in (21)) is convex and continuous. Unfortunately this result is not sufficient to justify that the problem (21) is well-posed. The difficulty comes from the fact that, even under Assumption 1, one can not guarantee that the problem

$$\text{Find } u_h \in V_h^k(\Omega) \text{ such that } \int_\Omega \varepsilon \nabla u_h \cdot \nabla v_h \, \mathrm{d}\boldsymbol{x} = \int_\Omega f v_h \, \mathrm{d}\boldsymbol{x} \quad \forall v_h \in V_h^k(\Omega)$$

is well-posed even for $h$ small enough. To cope with this difficulty, an idea is to use again the Tikhonov regularization (see §4.3), with a regularization parameter that depends now on $h$. This idea was originally proposed in [22] for the case of elliptic equations and then, was used by Assyr Abdulle et al. in [1] for the case of problems with sign-changing coefficients. Here, we explain how to adapt it to our approach. The idea is to replace the cost functional $J_{0,h}$ in (21) by the functional $J_h : V_h^k(\Omega_2) \to \mathbb{R}_+$ such that for all $w_h \in V_h^k(\Omega_2)$, we have

$$J_h(w_h) := \frac{1}{2} \int_\Sigma |u_h^{w_h} - u_{2,h}^{w_h}|^2 \, \mathrm{d}\sigma + \lambda_h \|w_h\|_{\tilde{\varepsilon}_1}^2,$$

where $\lambda_h > 0$ tends to zero as $h$ goes to 0. Since $\lambda_h > 0$ for all $h > 0$, the functional $J_h$ is strictly convex and coercive. This guarantees that the optimization problem

$$\min_{w_h \in V_h^k(\Omega_2)} J_h(w_h) \tag{22}$$

has a unique solution that we denote by $w_{k,h}^*$. All the difficulty now is to choose the parameter $\lambda_h$ in order to be able to ensure the convergence of $(w_{k,h}^*)_h$ towards a solution to (8) as $h$ tends to zero. This is the main goal of the next paragraph.

## 5.3 Convergence of the method

The starting point of our discussion is the following

**Lemma 5.1.** *We have the estimate*

$$J_h(w^*_{k,h}) \leq \frac{1}{2} \int_\Sigma |u_h^{w^*_J} - u_{2,h}^{w^*_J}|^2 \, d\sigma + \lambda_h \|w^*_J\|^2_{\tilde{\varepsilon}_1} \tag{23}$$

*where $w^*_J$ is defined by* (10).

*Proof.* Starting from the fact that $\pi_h^k w^*_J \in V_h^k(\Omega_2)$ and using that $w^*_{k,h}$ is the unique solution to the optimization problem (22), we conclude that $J_h(w^*_{k,h}) \leq J_h(\pi_h^k w^*_J)$. On the other hand, the identity (20) allows us to write

$$J_h(\pi_h^k w^*_J) = \frac{1}{2} \int_\Sigma |u_h^{w^*_J} - u_{2,h}^{w^*_J}|^2 \, d\sigma + \lambda_h \|\pi_h^k w^*_J\|^2_{\tilde{\varepsilon}_1}.$$

The Lemma is then proved by recalling the estimate (19). □

In order to simplify notations, for $h > 0$ and $w \in V_2(\Omega_2)$, we denote by $A_h(w)$ the real number

$$A_h(w) = \frac{1}{2} \int_\Sigma |u_h^w - u_{2,h}^w|^2 \, d\sigma.$$

From (20), we know that for all $w \in V_2(\Omega_2)$, we have $A_h(w) = J_0^h(\pi_h^k w)$. The main result of this paragraph is the following theorem.

**Theorem 5.1.** *Assume that the parameter $\lambda_h$ can be chosen such that the sequences $(\lambda_h)_h$ and $(A_h(w^*_J)/\lambda_h)_h$ tend to zero as $h$ tends to zero. Then, as $h$ goes to $0$:*

- *the sequence $(w^*_{k,h})_h$ converges to $w^*_J$ in $V_2(\Omega_2)$;*

- *the sequence $(u_h^{w^*_{k,h}})_h$ converges to $E_H(u_1)$ in $H_0^1(\Omega)$, resp. the sequence $(u_{2,h}^{w^*_{k,h}})_h$ converges to $u_2$ in $V_2(\Omega_2)$, where $(u_1, u_2)$ is the solution to (4) and $E_H(u_1)$ is the extension of $u_1$ defined in (12).*

*Proof.* The strategy of proof is similar to the one of proposition 4.3. To simplify notations, we denote by $u^{k,h} \in V_h^k(\Omega)$ and $u_2^{k,h} \in V_h^k(\Omega_2)$ the functions

$$u^{k,h} = u_h^{w^*_{k,h}} \qquad \text{and} \qquad u_2^{k,h} = u_{2,h}^{w^*_{k,h}}.$$

In order to make the proof as clear as possible, we divide it into four steps.

**Step 1: weak convergence of $(w^*_{k,h})_h$, $(u^{k,h})_h$ and $(u_2^{k,h})_h$.** Starting from the estimate

$$\|w^*_{k,h}\|^2_{\tilde{\varepsilon}_1} \leq J_h(w^*_{k,h})/\lambda_h \leq A_h(w^*_J)/\lambda_h + \|w^*_J\|^2_{\tilde{\varepsilon}_1} \tag{24}$$

and using the fact that $(A_h(w^*_J)/\lambda_h)_h$ tends to 0 as $h$ goes to 0, we infer that $(w^*_{k,h})_h$ is bounded in $V_2(\Omega_2)$. This implies that, up to a sub-sequence, $(w^*_{k,h})_h$ converges weakly to some $w_0 \in V_2(\Omega)$. For the reader's convenience, this sub-sequence is still denoted by $(w^*_{k,h})_h$.

Since the problems in (18) are uniformly elliptic with respect to $h$, we know that the sequence $(u^{k,h})_h$ (resp. $(u_2^{k,h})_h$) converges weakly in $H_0^1(\Omega)$ (resp. in $V_2(\Omega_2)$) to some $u \in H_0^1(\Omega)$ (resp. $u_2 \in V_2(\Omega_2)$). Using the basic approximability property (17), we infer that $u = u^{w_0}$ and $u_2 = u_2^{w_0}$.

**Step 2: $w_0$ is a minimizer of $J$.** The continuity of trace operator and the compactness of the embedding $H^{1/2}(\Sigma) \subset L^2(\Sigma)$ ensure that

$$u^{k,h}_{|\Sigma} - u_2^{k,h}{}_{|\Sigma} \to u^{w_0}{}_{|\Sigma} - u_2^{w_0}{}_{|\Sigma}$$

in $L^2(\Sigma)$ as $h \to 0$. By noticing that

$$\frac{1}{2}\int_\Sigma |u^{k,h} - u_2^{k,h}|^2 \, d\sigma = J_0^h(w_{k,h}^*) \leq J_h(w_{k,h}^*) \leq \lambda_h(A_h(w_J^*)/\lambda_h + \|w_J^*\|_{\tilde\varepsilon_1}^2)$$

and using that $\lambda_h, A_h(w_J^*)/\lambda_h \to 0$ as $h$ goes to zero, we deduce that $u^{w_0}_{|\Sigma} - u_2^{w_0}{}_{|\Sigma} = 0$. This shows that $w_0$ is a minimizer of $J$.

**Step 3: strong convergence of $(w_{k,h}^*)_h$ to $w_J^*$.** Thanks to the fact that $A_h(w_J^*)/\lambda_h \to 0$ as $h \to 0$ and by means of the estimate (24), we can write

$$\limsup_{h \to 0} \|w_{k,h}^*\|_{\tilde\varepsilon_1} \leq \|w_J^*\|_{\tilde\varepsilon_1}.$$

On the other hand, since $(w_{k,h}^*)_h$ converges weakly to $w_0$ as $h \to 0$, we infer that

$$\|w_0\|_{\tilde\varepsilon_1} \leq \liminf_{h \to 0} \|w_{k,h}^*\|_{\tilde\varepsilon_1}.$$

This implies that $\|w_0\|_{\tilde\varepsilon_1} \leq \|w_J^*\|_{\tilde\varepsilon_1}$. Since $w_0$ is a minimizer of $J$, we conclude that $w_0 = w_J^*$. Furthermore, we also deduce that

$$\lim_{h \to 0} \|w_{k,h}^*\|_{\tilde\varepsilon_1} = \|w_0\|_{\tilde\varepsilon_1}.$$

As a result, by applying [9, Proposition III.32], we infer that $(w_{k,h}^*)_h$ converges, strongly, in $V_2(\Omega_2)$ to $w_0 = w_J^*$.

**Step 4: strong convergence of $(u^{k,h})_h$ and $(u_2^{k,h})_h$.** The ellipticity of the problems in (6), combined with the strong convergence of $(w_{k,h}^*)_h$ to $w_J^*$, imply the convergence of $(u^{k,h})_h$ in $H_0^1(\Omega)$ to $u^{w_J^*}$ and of $(u_2^{k,h})_h$ in $V_2(\Omega_2)$ to $u_2^{w_J^*}$.

The Lemma is then proved by using that $u^{w_J^*} = E_H(u_1)$ and by observing that these limits are independent of the chosen sub-sequences. $\square$

The rest of this paragraph is devoted to explain why it is possible to choose the parameter $\lambda_h$ in such a way that $(\lambda_h)_h$ and $(A_h(w_J^*)/\lambda_h)_h$ both converge to 0 as $h$ tends to 0. To do so, one needs to study the behaviour of $A_h(w_J^*)$ as $h$ tends to 0. Let us start with the following

**Proposition 5.1.** *Suppose that the coefficients $\tilde\varepsilon_1$ and $\varepsilon_2$ are smooth, or piecewise smooth. Assume that the solution $u$ to (1) belongs to $PH^{1+s}(\Omega)$ for some $s > 0$. Then there exists $s' \in (0, s]$ that depends only on the geometry of $\Omega_2$ and on the coefficient $\varepsilon_2$, and there exists $\sigma \in (0, 1]$ that depends only on the geometry of $\Omega$ and of $\Omega_2$ such that*

$$\|u^{w_J^*} - u_h^{w_J^*}\|_{H_0^1(\Omega)} \leq Ch^{p'}\|u\|_{PH^{1+p'}(\Omega)} \quad \text{and} \quad \|u_2^{w_J^*} - u_{2,h}^{w_J^*}\|_{\tilde\varepsilon_1} \leq Ch^{p'}\|u_2\|_{H^{1+p'}(\Omega_2)},$$

$$\|u^{w_J^*} - u_h^{w_J^*}\|_{L^2(\Omega)} \leq Ch^{p'+\sigma}\|u\|_{H^{1+p'}(\Omega)} \quad \text{and} \quad \|u_2^{w_J^*} - u_{2,h}^{w_J^*}\|_{L^2(\Omega_2)} \leq Ch^{p'+\sigma}\|u_2\|_{H^{1+p'}(\Omega_2)}$$

*with $C$ independent of $h$ and $p' = \min(s', k)$.*

*Proof.* Along this proof, $C$ denotes a positive constant whose value can change from one line to the next but does not depend on $h$.

Given that $u^{w_J^*} = E_H(u_1)$ solves (12), and since $u_1 \in H^{1+s}(\Omega_1)$, it follows that $E_H(u_1)_{|\Omega_2}$ exhibits some extra-regularity because $\varepsilon_2$ is (piecewise) smooth (via an ad hoc shift theorem). In other words, there exists $s' \in (0, s]$ such that $u^{w_J^*} \in PH^{1+s'}(\Omega)$.

Given that $u_2^{w_J^*} = u_2 \in PH^{1+s}(\Omega_2) \subset PH^{1+s'}(\Omega_2)$ and since the problems in (6) are elliptic with (piecewise) smooth coefficients $\tilde\varepsilon_1$ and $\varepsilon_2$, we obtain the estimates (see [15])

$$\|u^{w_J^*} - u_h^{w_J^*}\|_{H_0^1(\Omega)} \leq Ch^{p'}\|u\|_{PH^{1+p'}(\Omega)} \text{ and } \|u_2^{w_J^*} - u_{2,h}^{w_J^*}\|_{\tilde\varepsilon_1} \leq Ch^{p'}\|u_2\|_{H^{1+p'}(\Omega_2)},$$

where $p' = \min(s', k)$. By applying the classical Aubin-Nitsche's lemma (see [15, Theorem 3.2.4]), we infer that there exists $\sigma \in (0, 1]$ such that

$$\|u^{w_J^*} - u_h^{w_J^*}\|_{L^2(\Omega)} \leq Ch^{p'+\sigma}\|u\|_{H^{1+p'}(\Omega)} \text{ and } \|u_2^{w_J^*} - u_{2,h}^{w_J^*}\|_{L^2(\Omega_2)} \leq Ch^{p'+\sigma}\|u_2\|_{H^{1+p'}(\Omega_2)}.$$

$\square$

Now we have all the tools to study the behavior $A_h(w_J^*)$ as $h$ goes to 0.

**Corollary 5.1.** *Under the same assumptions as in proposition 5.1, one has*

$$A_h(w_J^*) \leq Ch^{2p'+\sigma}$$

*with C independent of h and $p' = \min(s', k)$.*

*Proof.* Applying the multiplicative trace inequality (recalled in proposition A.1) and using the estimates of proposition 5.1 yield the estimates

$$\|u^{w_J^*} - u_h^{w_J^*}\|_{L^2(\Sigma)}^2 \leq Ch^{2p'+\sigma}\|u\|_{PH^{1+p'}(\Omega)} \text{ and } \|u_2^{w_J^*} - u_{2,h}^{w_J^*}\|_{L^2(\Sigma)}^2 \leq Ch^{2p'+\sigma}\|u_2\|_{H^{1+p'}(\Omega_2)}.$$

By design, one has $u^{w_J^*}_{|\Sigma} = u_2^{w_J^*}_{|\Sigma}$. So, observing that

$$\|u_h^{w_J^*} - u_{2,h}^{w_J^*}\|_{L^2(\Sigma)}^2 \leq 2(\|u^{w_J^*} - u_h^{w_J^*}\|_{L^2(\Sigma)}^2 + \|u_2^{w_J^*} - u_{2,h}^{w_J^*}\|_{L^2(\Sigma)}^2),$$

we conclude that $A_h(w_J^*) \leq Ch^{2p'+\sigma}$. $\square$

The previous result gives us a simple way to choose the parameter $\lambda_h$ in order to ensure that both $(\lambda_h)_h$ and $(A_h(w_J^*)/\lambda_h)_h$ tend to 0 as $h$ tends to 0.

**Corollary 5.2.** *Under the same assumptions as in proposition 5.1, any parameter $\lambda_h$ of the form $\lambda_h = Ch^q$ with $C > 0$ independent of h and $q \in (0, 2p' + \sigma)$ satisfies the conditions of theorem 5.1.*

**Remark 5.2.** *It is worth to note that the value of $s'$ is prescribed only by the regularity of $E_H(u_1)$, the harmonic-like extension $u_1$ that satisfies (12). Assume for instance that $\Omega$ and $\Omega_2$ are convex domains, and that the coefficients $\varepsilon_1$ and $\varepsilon_2$ are constant. In this case, one can choose $\tilde{\varepsilon}_1$ to be constant over $\Omega$. Let us consider the case where the solution u to (1) belongs to $PH^{1+s}(\Omega)$ for any $s \in (0, \sigma_D(\varepsilon))$. Then, because $E_H(u_1)_{|\Omega_2} \in H^1(\Omega_2)$ is governed by (cf. (12)): $-\Delta(E_H(u_1)) = 0$ in the convex domain $\Omega_2$ with Dirichlet data in $H^{1/2+s}(\partial\Omega_2)$, one has $E_H(u_1)_{|\Omega_2} \in H^{1+s}(\Omega_2)$. In other words, $s' = s$, and $p' = \min(s', k) = s' = s$. Also, because $\Omega$ and $\Omega_2$ are convex, one finds that $\sigma = 1$. Hence, in the statement of corollary 5.2, one may choose any $q \in (0, 2\sigma_D(\varepsilon) + 1)$.*

Thanks to theorem 5.1, using the conditions of corollary 5.2, one obtains the convergence of the discrete solutions to the exact solution.
On the one hand, convergence is guaranteed even on meshes that are not T-conforming. And, compared to [1], convergence holds in very general situations, namely as soon as there is a shift theorem for problem (1), cf. theorem 2.1, even with a regularity exponent $\sigma_D(\varepsilon) < 1/2$.
On the other hand, there is no associated convergence rate. Assuming a Céa lemma-like result, and using the same notations as above, the *expected* convergence rate is $h^{p'}$ in $H_0^1$-norm, and $h^{p'+\sigma}$ in $L^2$-norm. Whereas, classically, the *optimal* convergence rate is $h^k$ in $H_0^1$-norm, and $h^{k+1}$ in $L^2$-norm.

## 6 Numerical experiments

In this section we turn our attention to the validation of the numerical method that we have proposed. We limit ourselves to the case of 2D domains and use $P^1$ Lagrange finite elements. The numerical results that we present below have been obtained with the help of the library `FreeFem++`[2]. In particular, to solve the optimization problem (22), we used the `BFGS` function. In all the numerical experiments presented below, we have used the `BFGS` function with the following parameters: $eps = 10^{-6}, nbiter = 10, nbiterline = 1$ (see the `FreeFem++` documentation for more details).

---
[2] See https://freefem.org/.

Since the well-posedness of (1) depends on the shape of the interface $\Sigma$, we test the performance of our method in three different configurations. In the first one, $\Sigma$ is flat, in the second one, $\Sigma$ is circular interface and in the last one, $\Sigma$ has a "corner", in the sense that the angle at the intersection with the boundary is not a right angle. In all these experiments, we suppose that the coefficients $\varepsilon_1$ and $\varepsilon_2$ are constant with $\varepsilon_1 = 1$. We denote by $\kappa_\varepsilon$ the contrast $\kappa_\varepsilon = \varepsilon_2/\varepsilon_1$.

The shape, smoothness and (respective) volumes of $\Omega_1$ and $\Omega_2$ are taken into account to choose the domain $\Omega^\star \in \{\Omega_1, \Omega_2\}$ to which the extension is performed (we recall that one must have $\text{meas}_{\partial\Omega}(\partial\Omega^\star \setminus \Sigma) > 0$, see footnote[1] on page 4). Indeed, to have a better convergence rate, one should choose $\Omega^\star$ convex, or with as smooth a boundary as possible. Also, in order to speed up the convergence of the `BFGS` function, we must choose $\Omega^\star$ as small as possible. Once $\Omega^\star$ is fixed, one has to extend the function $\varepsilon_1$ or $\varepsilon_2$ to all the domain $\Omega$. Because the coefficients are constant, we extend $\varepsilon_1$ (resp. $\varepsilon_2$) by $\varepsilon_1$ in $\Omega_2$ (resp. in $\Omega_1$). In the case where $\Sigma$ is flat or circular, we take $\Omega^\star = \Omega_2$. In the third configuration, we take $\Omega^\star = \Omega_1$.

### 6.1 Flat interface

In this paragraph, we take

$$\Omega_1 = \{(x,y) \in (0; 1/2) \times (0; 1)\} \quad \text{and} \quad \Omega_2 = \{(x,y) \in (1/2; 1) \times (0; 1)\}$$

(a flat interface and a domain which is symmetric with respect to $\Sigma$). We consider a family of meshes of $\overline{\Omega}$ satisfying Assumption 2 (see Figure 2). In the rest of this paragraph we suppose that $\kappa_\varepsilon \neq -1$. To test the performance of our method, we work with the same example considered in [1, 14]. Define the function $u_{\kappa_\varepsilon}$ such that

$$u_{\kappa_\varepsilon}(x,y) = \begin{cases} (x^2 + bx)\sin(\pi y) & \text{if } x < 1/2 \\ a(x-1)\sin(\pi y) & \text{if } 1/2 < x \end{cases}, \quad \text{where } a = \frac{1}{2(\kappa_\varepsilon + 1)} \text{ and } b = -\frac{\kappa_\varepsilon + 2}{2(\kappa_\varepsilon + 1)}.$$

and consider it as an exact solution to (1). This is possible because $\text{div}(\varepsilon \nabla u_{\kappa_\varepsilon}) \in L^2(\Omega)$. The source term $f$ is computed accordingly. As observed in remark 5.2, by choosing $\lambda_h = Ch^q$ with $q \in (0,3)$, the method is convergent. In our experiment, we take $\lambda_h = 0.002h^2$. We work with $\kappa_\varepsilon = -1.001$. The behaviors of the relative $L^2$-norm error ($\|e_h^r\|_0$) and the relative $H_0^1$-norm error ($\|e_h^r\|_1$) between the exact solution and the numerical one are reported in Figure 2. We observe that both rates of convergence are equal to 2.

**Remark 6.1.** *The constant $C$ in $\lambda_h = Ch^q$ must be adjusted by the user according to the contrast $\kappa_\varepsilon$ in order to obtain a fast convergence of the method. Clearly this depends on $\|w_J^*\|_{\tilde{\varepsilon}_1}$. Using the fact that $w_J^* = w_H$ and owing to (15) we see that this depends on the jump of the normal derivative (across $\Sigma$) between $u_1$ and its harmonic extension. It is also important to note that, once $q$ is fixed and when $h$ is small enough, the choice of $C$ does not affect the convergence of the method.*

### 6.2 The case of a circular interface

In this paragraph, we consider the case where the domains $\Omega_1$ and $\Omega_2$ are such that $\Omega_1 = \{\boldsymbol{x} \in \mathbb{R}^2 \mid |\boldsymbol{x}| < 1\}$ and $\Omega_2 = \{\boldsymbol{x} \in \mathbb{R}^2 \mid 1 < |\boldsymbol{x}| < 2\}$. In proposition A.2, we prove that $A_\varepsilon$ is an isomorphism $\kappa_\varepsilon \notin \{-1\} \cup \mathscr{S}$ with $\mathscr{S} := \{-(1-(1/2)^{2n})/(1+(1/2)^{2n}) \mid n \in \mathbb{N}^*\}$. For this reason, we consider the case where $\kappa_\varepsilon = -2 \notin \mathscr{S}$. Given that both $\Omega_2$ and $\Omega$ have smooth boundaries, we infer that $\sigma = 1$ and $s' = s$. By taking $f$ as the source term associated to the function

$$u_{\kappa_\varepsilon}(x,y) = \begin{cases} r^2 + b & \text{if } r < 1 \\ a(r-2)^2 & \text{if } 1 < r < 2. \end{cases}, \quad \text{with } r = \sqrt{x^2 + y^2}, a = -1/\kappa_\varepsilon \text{ and } b = a - 1$$

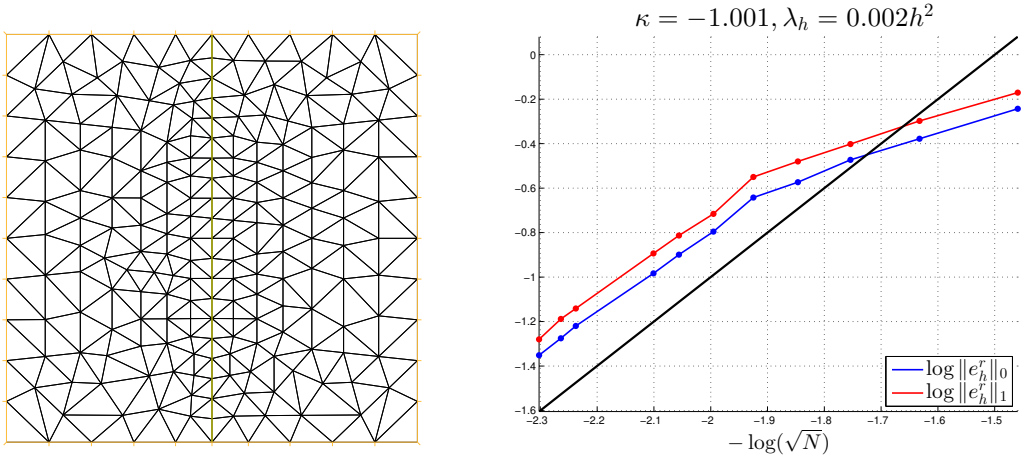

Figure 2: A given mesh (left). Behavior of the relative $L^2$ and $H_0^1$ errors with respect to the meshsize $h \sim \sqrt{N}$, where $N$ is the total number of nodes of the mesh (right).

and by taking $\lambda_h = 0.002h^2$, we obtain the results displayed in Figure 3. We observe that the method converges with optimal rate (ie. the relative $L^2$-norm error ($\|e_h^r\|_0$) is of order 2, while the relative $H_0^1$-norm error is of order 1), even though the exterior boundary and the interface are curved.

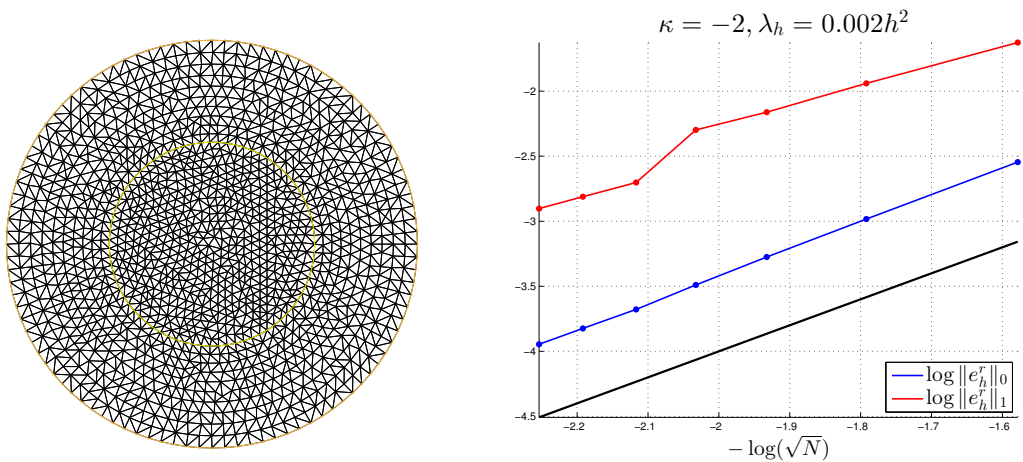

Figure 3: A given mesh (left). Behavior of the relative $L^2$ and $H_0^1$ errors with respect to the meshsize $h \sim \sqrt{N}$, where $N$ is the total number of nodes of the mesh (right)

### 6.3 The case of an interface with corner

Now, we consider the configuration where the interface $\Sigma$ has a corner. More precisely, we assume that $\Omega := \{\boldsymbol{x} \in \mathbb{R}^2 \mid |\boldsymbol{x}| < 1 \text{ and } \arg(\boldsymbol{x}) \in (0; \pi/2)\}$ and $\Omega_1 := \{\boldsymbol{x} \in \Omega \mid \arg(\boldsymbol{x}) \in (0; \pi/4)\}$ (see Figure 4). In such configuration, it can be proved (see [4]) that $A_\varepsilon$ is an isomorphism if and only if $\kappa_\varepsilon \in \mathbb{R}_-^* \setminus [-3, -1]$. Furthermore, contrarily to the two previous cases, in this configuration the solution to (1) can be very singular near the origin. Indeed, it was proved in [12, Chapter 2] that the regularity of the solution to (1) depends in $\kappa_\varepsilon$ and can be very low as $\kappa_\varepsilon$ approaches $[-3, -1]$: more precisely,

$$\lim_{\kappa_\varepsilon \to -3^-} \sigma_D(\varepsilon) = \lim_{\kappa_\varepsilon \to -1^+} \sigma_D(\varepsilon) = 0.$$

As a matter of fact, the value of the regularity exponent $\sigma_D(\varepsilon)$ is $\Re e(\lambda_0)$, where $\lambda_0$ is the solution to

$$\kappa_\varepsilon = -\tan(3\lambda\pi/4)/\tan(\lambda\pi/4) \tag{25}$$

that has the smallest positive real part. Note that one can show (see [12, Chapter 3]) that all the solutions to (25) are real-valued. In the particular cases where $\kappa_\varepsilon = -5$ and $\kappa_\varepsilon = -3.1$, one finds, respectively that $\lambda_0 \approx 0.458$ and $\lambda_0 \approx 0.139$. As mentioned previously this regularity result is optimal. Indeed, one can check that the function

$$u_{\lambda_0}(r,\theta) := (1-r)r^{\lambda_0} \begin{cases} \sin(\lambda_0 \theta)/\sin(\lambda_0 \pi/4) & \theta \in (0; \pi/4), \\ \sin(\lambda_0(\pi-\theta))/\sin(3\lambda_0\pi/4) & \theta \in (\pi/4; \pi) \end{cases}$$

satisfies $\text{div}(\varepsilon \nabla u_{\lambda_0}) \in L^2(\Omega)$. Observe that $u_{\lambda_0} \notin \text{PH}^{\lambda_0}(\Omega)$. This means that $u_{\lambda_0} \notin \text{PH}^{3/2}(\Omega)$. Now, given that $\Omega$ and $\Omega_2$ are both convex, owing to proposition 5.1, we can say that by choosing $\lambda_h = Ch^q$ with $q \in (0, 3\lambda_0)$, the convergence of the method can be guaranteed. In the case $\kappa_\varepsilon = -5$ (resp. $\kappa_\varepsilon = -3.1$), we work with $\lambda_h = h^{1.3}$ (resp. $\lambda_h = h^{0.4}$).

The behaviors of the relative $L^2$-norm error and the relative $H_0^1$-norm error (for the cases $\kappa_\varepsilon = -5$ and $\kappa_\varepsilon = -3.1$) are given in Figure 4. In either case, the expected rate of convergence is equal to $\lambda_0$ ($\approx 0.458$ when $\kappa_\varepsilon = -5$ and $\approx 0.139$ when $\kappa_\varepsilon = -3.1$) for the case of the $H_0^1$-norm error, while it is equal to $2\lambda_0$ for the case of the $L^2$-norm error. We observe that, in both cases, the method converges with optimal rate of convergence.

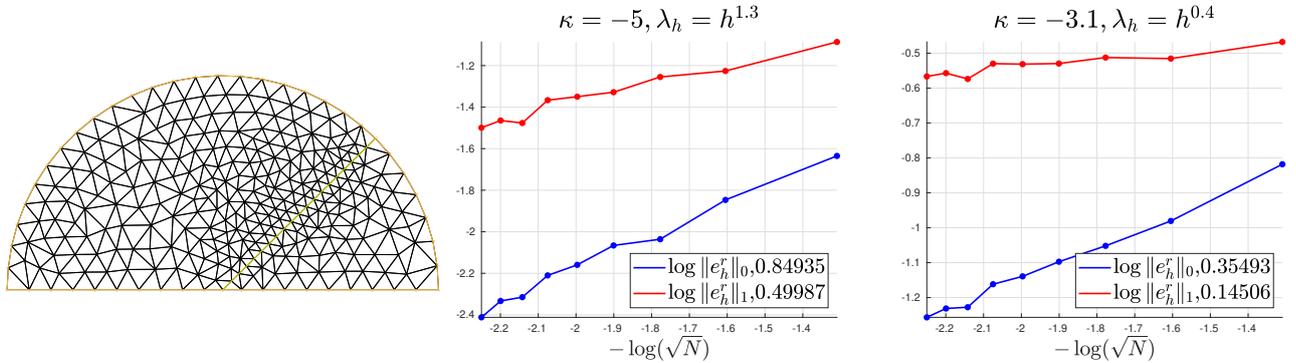

Figure 4: A given mesh (left). Behavior of the relative $L^2$ and $H_0^1$ errors with respect to the mesh-size $h \sim \sqrt{N}$, where $N$ is the total number of nodes of the mesh, with the observed convergence rates in the caption, when $\kappa_\varepsilon = -5$ (center) and $\kappa_\varepsilon = -3.1$ (right).

# 7 Concluding remarks

In this work, we have presented a new numerical method to approximate the solution to the scalar transmission problem with sign-changing coefficients. We proved that the method converges without any restriction on the mesh sequence, nor on the regularity of the solution. This result has been illustrated by several 2D numerical experiments. The convergence rate of our method seems to be optimal. In order to improve the performance of the method, several questions can be studied:

1. Choose the parameter $\lambda_h$ in order to accelerate convergence. An interesting idea would be to find an adaptive approach to fit its value. Also, one could use adaptive mesh refinement, together with a posteriori estimates. We refer to [17] for estimators that deliver guaranteed error bounds, and that are robust with respect to the sign-changing coefficient $\varepsilon$.

2. Work with other regularization approaches, ie. other choices for the coefficient $\tilde\varepsilon_1$, and/or an alternative to the Tikhonov regularization method.

3. In the case where the interface is not regular, it would be interesting to combine our approach with other existing methods for solving PDE with singular solution such as the Singular Complement Method [16].

Besides that, it will be also interesting to extend this approach to other models involving sign-changing coefficients.

# Acknowledgments

The authors would like to thank Anne-Sophie Bonnet-Ben Dhia and Lucas Chesnel for their helpful comments on the first draft of this article.

# A    Missing results

*Proof of Proposition 4.3 .* From the definition of $w_\delta^*$, we can write that

$$\delta \|w_\delta^*\|_{\tilde{\varepsilon}_1}^2 \leq J^\delta(w_\delta^*) \leq J^\delta(w_J^*) = J(w_J^*) + \delta \|w_J^*\|_{\tilde{\varepsilon}_1}^2 = \delta \|w_J^*\|_{\tilde{\varepsilon}_1}^2.$$

This means that for all $0 < \delta$, there holds $\|w_\delta^*\|_{\tilde{\varepsilon}_1} \leq \|w_J^*\|_{\tilde{\varepsilon}_1}$. As a result $(w_\delta^*)$ is bounded in $V_2(\Omega_2)$. This implies that, up to a sub-sequence, $(w_\delta^*)_\delta$ converges, as $\delta$ tends to 0, weakly in $V_2(\Omega_2)$ to some $w_0 \in V_2(\Omega_2)$. For the reader's convenience, this sequence is also denoted by $(w_\delta^*)_\delta$. Now, let us prove that $w_0$ is a minimizer of $J$. To do that, we start by observing that for all $\delta > 0$, we have

$$0 \leq J(w_\delta^*) \leq J^\delta(w_\delta^*) \leq J^\delta(w_J^*) = \delta \|w_J^*\|_{\tilde{\varepsilon}_1}^2.$$

This shows that $(J(w_\delta^*))_\delta$ converges to zero as $\delta$ tends to zero. On the other hand, by using the result of proposition 4.1, we know that $(J(w_\delta^*))_\delta$ converges to $J(w_0)$. Consequently, $J(w_0) = 0$ and then $w_0$ is a minimizer of $J$.

The next step is to show that the convergence of $(w_\delta^*)_\delta$ to $w_0$ occurs in the strong sense and that $w_0 = w_J^*$. To do so, we observe that

$$\left| \begin{array}{l} \|w_\delta^*\|_{\tilde{\varepsilon}_1} \leq \|w_J^*\|_{\tilde{\varepsilon}_1} \ \forall \delta \quad \implies \quad \limsup_{\delta \to 0} \|w_\delta^*\|_{\tilde{\varepsilon}_1} \leq \|w_J^*\|_{\tilde{\varepsilon}_1} \\ \\ w_\delta^* \rightharpoonup w_0 \text{ in } V_2(\Omega_2) \quad \implies \quad \|w_0\|_{\tilde{\varepsilon}_1} \leq \liminf_{\delta \to 0} \|w_\delta^*\|_{\tilde{\varepsilon}_1}{}^3. \end{array} \right.$$

This implies that $\|w_0\|_{\tilde{\varepsilon}_1} \leq \|w_J^*\|_{\tilde{\varepsilon}_1}$. Thanks to the definition of $w_J^*$, we deduce that $w_0 = w_J^*$. With this in mind and with the help of the previous inequality, we conclude that

$$\lim_{\delta \to 0} \|w_\delta^*\|_{\tilde{\varepsilon}_1} = \|w_J^*\|_{\tilde{\varepsilon}_1}.$$

Since $V_2(\Omega_2)$ is a Hilbert space, it follows (see [9, Proposition III.32]) that $w_\delta \to w_J^*$ in $V_2(\Omega_2)$. By noticing that $w_J^*$ is independent of the considered sub-sequence, the result is then proved. □

**Proposition A.1.** *[8, Theorem 1.6.6] Let $\Omega$ be an open, bounded, connected subset of $\mathbb{R}^d$ ($d = 2, 3$) with a Lipschitz boundary. Then the estimate*

$$\|u\|_{L^2(\partial\Omega)} \leq C \|u\|_{L^2(\Omega)}^{1/2} \|u\|_{H^1(\Omega)}^{1/2} \quad \forall u \in H^1(\Omega)$$

*holds with $0 < C$ independent of $u$.*

**Proposition A.2.** *Let $\Omega_1 = \{\boldsymbol{x} \in \mathbb{R}^2 \mid |\boldsymbol{x}| < 1\}$ and $\Omega_2 = \{\boldsymbol{x} \in \mathbb{R}^2 \mid 1 < |\boldsymbol{x}| < 2\}$. Assume that $\kappa_\varepsilon := \varepsilon_2/\varepsilon_1 \notin \{-1\} \cup \mathscr{S}$ with*

$$\mathscr{S} := \left\{ -\frac{1 - (1/2)^{2n}}{1 + (1/2)^{2n}} \mid n \in \mathbb{N}^* \right\}.$$

*Then the operator $A_\varepsilon : H_0^1(\Omega) \to H_0^1(\Omega)$ is an isomorphism.*

---

[3]This is a consequence of the fact that the norm of a Banach space is weakly lower semicontinuous, see [9, Proposition III.5 (iii)].

**Remark A.1.** *Note that in accordance with the results concerning the Neumann-Poincaré operator [24, Chapter 1], we observe that $-1$ is an accumulation point of $\mathscr{S}$.*

*Proof.* [12, Theorem 1.3.3] guarantees that $A_\varepsilon$ is Fredholm of index 0 when $\kappa_\varepsilon \neq -1$. Therefore it suffices to study its kernel. Let $u \in \mathrm{H}_0^1(\Omega)$ be such that $A_\varepsilon u = 0$. Then $u_1 := u_{|\Omega_1}$ and $u_2 = u_{|\Omega_2}$ satisfy

$$\begin{cases} \Delta u_1 = 0 & \text{in } \Omega_1 \\ \Delta u_2 = 0 & \text{in } \Omega_2 \\ u_1(1,\theta) = u_2(1,\theta) & \text{and} \quad \partial_r u_1(1,\theta) = \kappa_\varepsilon \partial_r u_2(1,\theta) \quad \forall \theta \in [0; 2\pi]. \end{cases}$$

Since the problem is invariant with respect to $\theta$, by Fourier decomposition for $u_1$, $u_2$ we have the representations:

$$u_1(r,\theta) = \sum_{n \in \mathbb{N}} a_n r^n \mathrm{e}^{in\theta} \quad \text{and} \quad u_2(r,\theta) = b_0 \ln(r/2) + \sum_{n \in \mathbb{Z}^*} b_n((r/2)^n - (r/2)^{-n}) \mathrm{e}^{in\theta},$$

where $a_n, b_n \in \mathbb{C}$. Using the transmission conditions, we get

$$\begin{vmatrix} a_0 = b_0 \ln(1/2), & 0 = b_0 \kappa_\varepsilon \\ a_n = b_n((1/2)^n - (1/2)^{-n}), & a_n = b_n((1/2)^n + (1/2)^{-n})\kappa_\varepsilon, & n \in \mathbb{N}^* \\ 0 = b_n((1/2)^n - (1/2)^{-n}), & 0 = b_n((1/2)^n + (1/2)^{-n})\kappa_\varepsilon, & -n \in \mathbb{N}^*. \end{vmatrix}$$

Therefore we deduce that $A_\varepsilon$ is injective when $\kappa_\varepsilon \notin \mathscr{S}$. $\square$

## B  On the use of the adjoint approach to compute the gradient of the cost functional $J$

The adjoint approach was introduced in [11] as a method for computing the gradient of cost functions that depend in non-explicit way of the main variable of the problem, namely via the solution of PDEs (the state equations) in which the main variable plays the role of a parameter. Here, we are going to explain how to apply this method to our case. The idea is to introduce a Lagrangian functional $\mathscr{L} : \mathrm{V}_2(\Omega_2) \times \mathrm{H}_0^1(\Omega) \times \mathrm{V}_2(\Omega_2) \times \mathrm{H}_0^1(\Omega) \times \mathrm{V}_2(\Omega_2) \to \mathbb{R}$ such that

$$\mathscr{L}(w, u, u_2, g, g_2) = \frac{1}{2} \int_\Sigma |u - u_2|^2 \, \mathrm{d}\sigma + a_1(w, u, g) + a_2(w, u, u_2, g_2)$$

in which $a_1(w, u_1, g)$ and $a_2(w, u_2, g_2)$ are respectively given by

$$\begin{vmatrix} a_1(w, u, g) = \int_\Omega \tilde{\varepsilon}_1 \nabla u \cdot \nabla g \, \mathrm{d}\boldsymbol{x} - \int_{\Omega_1} fg \, \mathrm{d}\boldsymbol{x} - \int_{\Omega_2} \tilde{\varepsilon}_1 \nabla w \cdot \nabla g \, \mathrm{d}\boldsymbol{x} \\ a_2(w, u, u_2, g_2) = \int_{\Omega_2} \varepsilon_2 \nabla u_2 \cdot \nabla g_2 \, \mathrm{d}\boldsymbol{x} - \int_{\Omega_2} f_2 g_2 \, \mathrm{d}\boldsymbol{x} + \int_{\Omega_2} \tilde{\varepsilon}_1 \nabla (w - u) \cdot \nabla g_2 \, \mathrm{d}\boldsymbol{x}. \end{vmatrix}$$

The functions $g \in \mathrm{H}_0^1(\Omega), g_2 \in \mathrm{V}_2(\Omega_2)$ are the adjoint variables associated to $u, u_2$ respectively. Let $(u^w, u_2^w)$ be the solution to (6). By design, $a_1(w, u^w, g) = 0$ for all $g \in \mathrm{H}_0^1(\Omega)$, and $a_2(w, u^w, u_2^w, g_2) = 0$ for all $g_2 \in \mathrm{V}_2(\Omega_2)$, so one has

$$\mathscr{L}(w, u^w, u_2^w, g, g_2) = J(w) \qquad \forall g \in \mathrm{H}_0^1(\Omega), \ \forall g_2 \in \mathrm{V}_2(\Omega_2). \tag{26}$$

Clearly, the functional $\mathscr{L}$ is differentiable with respect to all its variables. For all $\mathbf{v} = (w, u, u_2, g, g_2) \in \mathrm{V}_2(\Omega_2) \times \mathrm{H}_0^1(\Omega) \times \mathrm{V}_2(\Omega_2) \times \mathrm{H}_0^1(\Omega) \times \mathrm{V}_2(\Omega_2)$, the partial derivatives of $\mathscr{L}$ at $\mathbf{v}$ belong respectively to $\partial_w \mathscr{L}(\mathbf{v}) \in (\mathrm{V}_2(\Omega_2))^*$, $\partial_u \mathscr{L}(\mathbf{v}) \in (\mathrm{H}_0^1(\Omega))^*$, $\partial_{u_2} \mathscr{L}(\mathbf{v}) \in (\mathrm{V}_2(\Omega_2))^*$, $\partial_g \mathscr{L}(\mathbf{v}) \in (\mathrm{H}_0^1(\Omega))^*$, $\partial_{g_2} \mathscr{L}(\mathbf{v}) \in (\mathrm{V}_2(\Omega_2))^*$.

Let $g \in H_0^1(\Omega)$ and $g_2 \in V_2(\Omega_2)$ be given, and $\mathbf{v}^w = (w, u^w, u_2^w, g, g_2)$. By taking the derivative of the relation (26) with respect to $w$, we find that, by applying the chain rule formula,

$$(J'(w), h)_{\tilde{\varepsilon}_1} = \langle \partial_w \mathscr{L}(\mathbf{v}^w), h \rangle + \langle \partial_u \mathscr{L}(\mathbf{v}^w), \frac{du^w}{dw}(h) \rangle + \langle \partial_{u_2} \mathscr{L}(\mathbf{v}^w), \frac{du_2^w}{dw}(h) \rangle, \quad \forall h \in V_2(\Omega_2).$$

Now, if there exists $(g^w, g_2^w) \in H_0^1(\Omega) \times V_2(\Omega_2)$ for which the equations

$$\partial_u \mathscr{L}(w, u^w, u_2^w, g^w, g_2^w) = 0 \quad \text{and} \quad \partial_{u_2} \mathscr{L}(w, u^w, u_2^w, g^w, g_2^w) = 0$$

are satisfied for all $w \in V_2(\Omega_2)$, this yields

$$(J'(w), h)_{\tilde{\varepsilon}_1} = \langle \partial_w \mathscr{L}(w, u^w, u_2^w, g_1^w, g_2^w), h \rangle \quad \forall w \in V_2(\Omega_2), \ \forall h \in V_2(\Omega_2).$$

To investigate the existence of $(g^w, g_2^w)$, we need to write down the expression of $\partial_u \mathscr{L}(\mathbf{v}^w)$ and $\partial_{u_2} \mathscr{L}(\mathbf{v}^w)$: By a direct calculus, one checks that

$$\langle \partial_u \mathscr{L}(\mathbf{v}^w), v \rangle = \int_\Omega \tilde{\varepsilon}_1 \nabla g \cdot \nabla v \, d\boldsymbol{x} - \int_{\Omega_2} \tilde{\varepsilon}_1 \nabla g_2 \cdot \nabla v \, d\boldsymbol{x} + \int_\Sigma (u^w - u_2^w) v \, d\sigma, \quad \forall v \in H_0^1(\Omega)$$

$$\langle \partial_{u_2} \mathscr{L}(\mathbf{v}^w), v_2 \rangle = \int_{\Omega_2} \varepsilon_2 \nabla g_2 \cdot \nabla v_2 \, d\boldsymbol{x} - \int_\Sigma (u^w - u_2^w) v_2 \, d\sigma \quad \forall v_2 \in V_2(\Omega_2).$$

Hence, the functions $(g^w, g_2^w) \in H_0^1(\Omega) \times V_2(\Omega_2)$ are governed by the following system of equations:

$$\left| \begin{array}{ll} \int_\Omega \tilde{\varepsilon}_1 \nabla g^w \cdot \nabla v \, d\boldsymbol{x} = \int_{\Omega_2} \tilde{\varepsilon}_1 \nabla g_2^w \cdot \nabla v \, d\boldsymbol{x} - \int_\Sigma (u^w - u_2^w) v \, d\sigma & \forall v \in H_0^1(\Omega) \\ \int_{\Omega_2} \varepsilon_2 \nabla g_2^w \cdot \nabla v_2 \, d\boldsymbol{x} = \int_\Sigma (u^w - u_2^w) v_2 \, d\sigma & \forall v_2 \in V_2(\Omega_2). \end{array} \right. \qquad (27)$$

Clearly the previous system of equations is well-posed. Therefore the functions $g^w$, $g_2^w$ are well-defined. We then have all the tools to prove the result stated in Lemma 4.1.

*Proof of Lemma 4.1.* Take $w \in V_2(\Omega_2)$. From the characterization (27) of $g^w$ and $g_2^w$, we deduce that for all $h \in V_2(\Omega_2)$, we have

$$(J'(w), h)_{\tilde{\varepsilon}_1} = \langle \partial_w \mathscr{L}(w, u^w, u_2^w, g_1^w, g_2^w), h \rangle.$$

On the other hand, one can compute explicitly the value of $\langle \partial_w \mathscr{L}(w, u, u_2, g, g_2), h \rangle$:

$$\langle \partial_w \mathscr{L}(w, u, u_2, g, g_2), h \rangle = \int_{\Omega_2} \tilde{\varepsilon}_1 \nabla h \cdot \nabla(g_2 - g_{|\Omega_2}) \, d\boldsymbol{x}.$$

This shows that $J'(w) = g^w{}_{|\Omega_2} - g_2^w$ and then the result is proved. $\square$